\documentclass[12pt]{article}

\usepackage[latin1]{inputenc}
\usepackage[T1]{fontenc}
\usepackage{times}

\usepackage{amssymb}
\usepackage{amsmath}
\usepackage{amscd}
\usepackage{epic}

\usepackage{xypic}

        \def\Z{\mathbb{Z}}

        \def\F{\mathbb{F}}

        \def\Afonts{\rmfamily \bfseries \normalsize}
        \def\afonts{\rmfamily \upshape \footnotesize}
        \def\Mclfonts{\rmfamily \bfseries \footnotesize}
        \def\mclfonts{\rmfamily \upshape \footnotesize}
        
        \newtheorem{Theo}{Theorem}[section]
	\newtheorem{Defi}[Theo]{Definition}
	\newtheorem{Propo}[Theo]{Proposition}
	\newtheorem{Coro}[Theo]{Corollary}
	\newtheorem{Lemme}[Theo]{Lemma}
	\newtheorem{ex}[Theo]{Example}
        \newtheorem{Rem}[Theo]{Remark}

\newcommand{\biindice}[3]%
 {
 
 {#1}_{\begin{array}[t]{c}
         {\scriptstyle #2} \\
         {\scriptstyle #3}
       \end{array}}
 
 }
\newcommand{\Deg }{\mathrm{Deg\, }}
\newcommand{\im }{\mathrm{Im\,}}
\newcommand{\Ker }{\mathrm{Ker\,}}
\newcommand{\norm }{\mathrm{norm\,}}
\newcommand{\tr }{\mathrm{trace\,}}
\newcommand{\Gal }{\mathrm{Gal\,}}
\newcommand{\ord }{\mathrm{ord\, }} 
\newenvironment{proof}{\noindent {{\bf Proof}.}}{\hspace*{\fill}$\square$\vskip 8pt}

         \def\bib#1#2#3#4{\bibitem{#1}{\it #2}{\ #3}{\ #4}}
\title{Descent of the Definition Field of a Tower of Function Fields and Applications}

\date{\today}
\author{St\'{e}phane Ballet$^*$, Dominique Le Brigand$^{**}$, Robert Rolland$^{***}$ }

\begin{document}
\maketitle


\bigskip

{\Afonts { \centerline {Abstract}}}
{\afonts Let us consider an algebraic function field
defined over a finite Galois extension $K$ of a perfect field
$k$. We give some conditions allowing
the descent of the definition field of the algebraic
function field  from $K$ to $k$. We apply these results to the descent
of the definition field of a tower of function fields.
We give explicitly the equations 
of the intermediate steps of an Artin-Schreier type extension reduced 
from $\F_{q^2}$ to $\F_q$.
By applying these results to a completed 
Garcia-Stichtenoth's tower
we improve the upper bounds and the upper asymptotic bounds of 
the bilinear complexity 
of the multiplication in finite fields.}


\medskip


\medskip

{\Mclfonts{{2000 Mathematics Subject Classification:}}} 
{\mclfonts {11G20, 14H25.}}


\renewcommand{\thefootnote}{\fnsymbol{footnote}}
\footnotetext{$^*$ S. Ballet, Laboratoire de G\'{e}om\'{e}trie Alg\'{e}brique 
et Applications \`{a} la
Th\'{e}orie de l'Information,
Universit\'{e} de la Polyn\'{e}sie Fran\c caise, BP 6570, 98702 Faa'a,
Tahiti, Polyn\'{e}sie Fran\c caise,\\
e-mail: ballet@upf.pf.}
\footnotetext{$^{**}$ D. Le Brigand, Institut de Math\'{e}matiques de Jussieu,
Universit\'{e} Pierre et Marie Curie, Paris VI, Campus Chevaleret,
175 rue du Chevaleret, F75013 Paris,\\
e-mail: Dominique.LeBrigand@math.jussieu.fr.}
\footnotetext{$^{***}$ R. Rolland, C.N.R.S. Institut de Math\'{e}matiques
de Luminy, case 930, F13288 Marseille cedex 9,\\
e-mail: rolland@iml.univ-mrs.fr.}

\section{Introduction}\label{un}
In the paper \cite{ballro} S. Ballet and R. Rolland have improved
the upper bounds and upper asymptotic bounds of the bilinear complexity 
of the multiplication in finite fields, for fields of characteristic $2$.
A key point, allowing the computation, is the existence of a descent of the 
definition field from $\F_{2^{2s}}$ to $\F_{2^{s}}$ 
for the Garcia-Stichtenoth's tower (cf. \cite{gast}) completed by
intermediate steps (cf. \cite{ball3}). The proof used in \cite{ballro}
cannot be performed in the odd characteristic case. 
So, it is of interest to know whether there exists, when the characteristic
is not $2$, a descent
of the definition field of a tower.

\medskip

In Section \ref{deux} we recall the notion of descent
of the definition field of an algebraic function field, and
we proceed with a brief study of this notion. We prove that
under certain conditions, the descent is possible. This part is closely
related to the descent theory, developed by A. Weil (cf. \cite{weil}).
However, because our goal is to study algebraic function fields
and not coordinate rings, we have replaced here the notion of descent
defined in \cite{serre} by a slightly weaker notion.

\medskip

In Section \ref{deubis}, we give explicitly the equations 
of the intermediate steps of an Artin-Schreier type extension reduced 
from $\F_{q^2}$ to $\F_q$.

\medskip

Section \ref{trois} is devoted
to the case of the completed
Garcia-Stichtenoth's tower. 
It contains a summary of the tower construction,
a proof that the descent of the definition field
from $F_{q^2}$ to $\F_q$ is possible. Using the results of Section  \ref{deubis}, 
we obtain an alternative proof
of the existence of the descent  and also we can give the  equation of each
step of the reduced tower. As an application
we derive an interesting bound for the bilinear complexity of the multiplication
in finite fields which extends the results of \cite{ballro}.

\section{General results on the descent of the definition field}\label{deux}

From now on, $k$ is a perfect field and $K$ is a finite Galois extension of $k$ such that 
$k\subset K\subset U$, where $U$ is a fixed algebraically closed field.
We will denote by $\Gamma=\Gal(K/k)$ the Galois group of this extension.
The group $\Gamma$ is also the Galois group of the Galois extension $K(x)/k(x)$
of the rational function field $k(x)$ by the rational
function field $K(x)$.
Let ${\cal C}$ be an absolutely irreducible algebraic curve.
We assume that ${\cal C}$
is defined over $K$. That is, the ideal of ${\cal C}$
can be generated by a family of polynomials 
which coefficients are in $K$.
Let $F=K({\cal C})$ be the 
algebraic function field of one variable associated to
${\cal C}$. We remark that 
$K$ is the full constant field of $F$.
The descent problem of the definition field of 
the curve ${\cal C}$ from $K$ to $k$ is the following: is it possible
to find an absolutely irreducible curve
${\cal C}'$ defined over $k$ and birationally isomorphic to ${\cal C}$?
If the answer is positive, let us set $F'=k({\cal C}')$. Then
$k$ is the full constant field of $F'$.
Algebraically speaking, the descent problem 
can be expressed in the following way:

\begin{Defi} Let $F/K$ be an algebraic function field of 
one variable with full constant field $K$.
We consider $F$ as a finite extension of $K(x)$ 
for some $x\in F$  transcendental over $K$. 
We say that the descent of the definition field of $F$ from $K$ to $k$
is possible,
if there exists an extension $F'$ of $k(x)$,
with full constant field $k$, such that
$F$ is isomorphic as $k$-algebra to $F'\otimes_k K$. We will say also
that we can reduce $F/K$ to $F'/k$.
\end{Defi}

For any algebraic function field $F/K$, we assume tacitely 
that $K$ is its full constant field and that $F$ is a subfield of
the algebraically closed field $U$.

\medskip

First we are going to study some simple general conditions
allowing the reduction. 
\begin{Theo}\label{desc} Let $F/K$ be an algebraic function field. 
If the extension $F/k(x)$ is Galois with Galois group ${\cal G}$ and 
if ${\cal G}$ is the semi-direct product ${\cal G}=G \rtimes \Gamma$,
where $G=\Gal(F/K(x))$ and $\Gamma=\Gal(K(x)/k(x))$,
then the descent of the definition field from
$K$ to $k$ is possible. 
\end{Theo}

\begin{proof}
Let us suppose that the extension $F/k(x)$ is Galois with
Galois group ${\cal G}$. Then $F/K(x)$ is Galois with
Galois group $G$ (cf. Figure (\ref{fig1})). As $K(x)$ is a Galois extension of $k(x)$,
the Galois group $G$ is a normal subgroup of ${\cal G}$ and $\Gamma$
is isomorphic to the factor group ${\cal G} / G$.

\begin{equation}
\label{fig1}
\xymatrix{
F' \ar@{--}[r] & F \\
k(x) \ar@{--}[u] \ar@{-}[r]_\Gamma \ar@{-}[ur]^{\cal G}&K(x)\ar@{-}[u]_G\\
k \ar@{-}[u] \ar@{-}[r]_\Gamma &K\ar@{-}[u]
}
\end{equation}

Then we have the exact sequence
$$1 \rightarrow G \overset{i}{\rightarrow} {\cal G} 
\overset{\underset{\pi}{\overset{s}{\curvearrowleft}}}{\rightarrow} \Gamma \rightarrow 1.$$
Moreover, ${\cal G}$ is the semi-direct product ${\cal G}=G \rtimes \Gamma$,
so that we can build a section $s$ lifting $\Gamma$ into ${\cal G}$.
We can say also that
for all $\sigma \in \Gamma=\Gal(K/k)$
there exists an automorphism $s({\sigma})$ of $F$ over $k(x)$
which is a continuation of $\sigma$. Moreover for each $\sigma, \tau$
in $\Gamma$, we have the condition
$$s(\sigma . \tau)=s(\sigma).s(\tau).$$
Hence, the group $\Gamma$ acts on $F$. 
More precisely its action $T$ is defined by
$$T(\sigma,x)=s(\sigma)(x).$$

Let $F'=F^{s(\Gamma)}$ be the fixed field 
of $s(\Gamma)$.
Since $K$ is a finite extension of $k$, $F$ is a finitely generated $k$-algebra
and any element of $F$ is integral over $F'$.
Then $F'$ is a finitely generated $k$-algebra (cf. \cite{serre}).
To prove that $F/K$ can be reduced,  
we have just to show that 
$F=F'\otimes_{k}K$. This is a consequence of Lemma 26 Chapter V
in \cite{serre}. 
\end{proof}

\begin{Coro} 
Suppose that $F/k(x)$ is an abelian extension with Galois group ${\cal G}$.  
If $\Gamma=\Gal(K/k)$ is a cyclic group of prime order $p$ and if 
$G=\Gal(F/K(x))$ is of order $m$ whith $\gcd(m,p)=1$ then 
the descent of the definition field from
$K$ to $k$ is possible.
\end{Coro}

\begin{proof}
 The commutative group ${\cal G}$ has $pm$ elements.
Hence it contains a cyclic subgroup $\Gamma'$ of order $p$. The morphism 
$s$ which sends a generator of $\Gamma$ on a generator of $\Gamma'$
defines a section. 
\end{proof}

\begin{Theo}\label{rdesc}
Let $F$ be a Galois extension of $K(x)$. Assume that 
it is possible to reduce $F/K$ to $F'/k$.
Then $F/k(x)$ is Galois (cf. Figure (\ref{fig1})). 
\end{Theo}

\begin{proof}
The extension $F/F'$ is Galois, and its Galois group is $\Gamma$. 
So that we can define an action from
$\Gamma$ on the Galois group $G=\Gal(F/K(x))$ by inner automorphisms:
$$\gamma(g)=\gamma^{-1}\circ g \circ \gamma.$$

Then, the semi-direct product $G\rtimes \Gamma$ is a group of 
$k(x)$-automorphisms 
of $F$. The cardinality of this group is exactly the degree of 
$F$ over $k(x)$. Hence we can conclude that $G\rtimes \Gamma$ is
the group of automorphisms of $F$ over $k(x)$ and that
the extension $F/k(x)$ is Galois, with Galois group ${\cal G}= G\rtimes \Gamma$.
\end{proof}


\begin{Theo}\label{princ}
Let $F/K$ and $E/K$ be algebraic function fields 
such that
$$K(x)  \subsetneq E \subsetneq  F.$$
We suppose that $F/K(x)$ is Galois. Let $H$
be the Galois group of 
$F/E$. 
Suppose that we can reduce $F/K$ 
to $F'/k$. 
Let $\Gamma$
be the Galois group of the extension $F/F'$.
Suppose that $\Gamma$ acts on $H$ by
inner automorphisms.
Then we can reduce $E/K$ to an algebraic function field $E'/k$.
\end{Theo} 

\begin{proof}
By Theorem \ref{rdesc} the extension $F/k(x)$ is Galois
with Galois group ${\cal G}$. 
Let ${\cal H}$ be the semi-direct product $H \rtimes \Gamma$ corresponding to
the action by inner automorphisms from $\Gamma$ on $H$.
This group is a subgroup of the Galois group ${\cal G}$. Hence there is
an algebraic function  field $E'/k$ corresponding to the field fixed by ${\cal H}$.
Note that $\Gamma$ is a subgroup of ${\cal H}$, so that
$E' \subset F'$. On the other hand the full constant field of $F'$ is $k$.
Hence, the full constant field of $E'$ is $k$.
The group $H$ is a subgroup of ${\cal H}$, then $H$ fixes $E'$
and consequently also $E'\otimes_kK$. We conclude that $E'\otimes_kK \subseteq E$.
But we have
$[E:K(x)]=\ord \left(G/H\right)$ and
$[E':k(x)]=\ord \left({\cal G}/{\cal H}\right).$
Hence, $[E:K(x)]=[E':k(x)]=[E'\otimes_kK:K(x)]$. 
Thus $E=E'\otimes_kK$, which completes the proof. 
\end{proof}
\begin{equation}
\label{fig3}
\xymatrix{
F^\prime \ar@{--}[r]^\Gamma
&F&F\ar@{-}[dd]^G \\
E' \ar@{--}[u]  \ar@{--}[r] \ar@{--}[ur]^{\cal H}&E\ar@{-}[u]_H&\\
k(x) \ar@{--}[u] \ar@{-}[uur]_<<<<<<<<<{\cal G} \ar@{-}[r]_\Gamma &K(x)\ar@{-}[u]&K(x)
}
\end{equation}

\begin{Rem}\label{rprinc} Let $E_1, E_2$ be two functions fields 
such that $K(x) \subset E_1 \subset E_2 \subset F$ and
satisfying the hypothesis of Theorem \ref{princ}. The construction
given in the proof of the previous theorem brings forth the descents $E'_1$
and $E'_2$ such that $k(x) \subset E'_1 \subset E'_2 \subset F'$.
\end{Rem}

\section{Explicit descent in particular cases}\label{deubis}

In this section, we assume that $q=p^n$ is a prime power.  

Our aim is the following. 
We consider an algebraic function field 
$L/\F_{q^2}$ 
and, 
for $u\in L$, we set $F=L(z)$, where $z^q+z=u$. We assume that 
there exists a place $\wp$ of $L$ such that 
$\nu_\wp(u)=-m$, with $m>0$ and $\gcd(m,q)=1$. 
Then the extension $F/L$ is of Artin-Schreier type (see \cite{gast}, Proposition 1.1) and, more precisely, $F/L$ is an elementary abelian extension of exponent $p$ and degree $q=p^n$.
Assume that $L/\F_{q^2}$ can be reduced over $\F_q$ and 
let $L'/\F_q$ is its reduced function field. 
Assume moreover that $u\in L'$ and set $G=L'(z)$, where $z^q+z=u$.
We want to prove that, for all $i$, $1\leq i<n$, there exists a 
subfield $G_i$
of $G$, such that $[G:G_i]=p^i$, and we want to give explicitly 
the equation of each function field $G_i/\F_q$.

Note that similar technics are used in \cite{stiGar} and \cite{Deo1}. 

\subsection{Linearized  or additive polynomials}


\begin{Defi} A linearized  or additive polynomial
$R(T)=\sum_{i=0}^{d}a_iT^{q^i}\in k[T]$  with coefficients in a field extension
$k$ of $\F_q$ is called a {\em $q$-polynomial over $k$}. If $R(T)\in\F_q[T]$ we say that $R(T)$ is a $q$-polynomial.
The {\em symbolic product} of two $q$-polynomials over $k$, $Q(T)$ and $M(T)$, is the $q$-polynomial over $k$ defined by
$$(M\star Q)(T)=M(Q(T)).$$
If $R(T)$ is a $q$-polynomial over $k$ such that $R(T)=(M\star Q)(T)$, 
where  $Q(T)$ and $M(T)$ are $q$-polynomials over $k$, we say that $M(T)$ 
divides symbolically $R(T)$.
\end{Defi}	
The symbolic product is  associative, distributive (with respect to the ordinary addition), but  it is not commutative if $\F_q\subsetneq k$ (see further).
To each $q$-polynomial over $k$
one can associate a $\F_q$-linear map from $k$ to itself.

Now, for $q$-polynomials, the symbolic product is  commutative and it is quite simple 
to know if a $q$-polynomial
divides symbolically another $q$-polynomial.
\begin{Lemme} Let $R(T)$ and $Q(T)$ be $q$-polynomials.
The following properties are equivalent:
\begin{itemize}
\item $Q(T)$ divides symbolically $R(T)$,
\item $Q(T)$ divides (in the ordinary sense) $R(T)$ in $\F_q[T]$.
\end{itemize}
\end{Lemme}
\begin{proof} see \cite[Theorem 3.62. p. 109]{Lidl}.
\end{proof}

\begin{ex}
\label{$Tq+T$} 
The polynomial $P(T)=T^q+T$ is a $q$-polynomial and also a $p$-polynomial. 
The set  of roots of $P$, ${\cal P}=\{\alpha\in\bar{\F}_q,\, P(\alpha)=0\}$ is a $\F_p$-vector space of dimension $n$ and
 $${\cal P}=\{\alpha\in \F_{q^2},\, \tr_{ \F_{q^2}/\F_q}(\alpha)=0\}.$$
If $p=2$, ${\cal P}=\F_q$ and otherwise ${\cal P}\subset \F_{q^2}$. Moreover ${\cal P}$ is stable under $\Gal(\F_{q^2}/\F_q)$, as well as any $\F_p$-subspace $H$ of $\cal P$, since $a^q=-a$ for all $a\in \cal P$.
Notice that, if $p=2$, the $\F_q$-linear map associated to $P(T)=T^q+T$ is the zero map from $\F_q$ to itself and that
the $p$-polynomials  $Q(T)=T^2+T$ and $M(T)=T^{2^{n-1}}+\cdots+T$ divide 
symbolically $P(T)$, since
$$P(T)=T^{2^n}+T=M(Q(T))=Q(M(T)).$$

\end{ex}

\begin{Defi} Let $k$ be a field extension of $\F_q$.
 We denote by $\tau_q\,: x\mapsto x^q$ the Frobenius endomorphism and we consider the $\F_q$-algebra   
$$k\{\tau_q\}=\left \{\sum_{i=0}^{d}a_i\tau_q^{i},\, a_i\in k\right\},$$
equipped with the composition law $\circ$.
\end{Defi}
Notice that, if $\alpha\in k$, we have  $\tau_q\circ(\alpha \tau_q^0)=\alpha^q\tau_q$, thus the composition law $\circ$ is commutative in $k\{\tau_q\}$ if and only if $k=\F_q$.
The $\F_q$-algebra of $q$-polynomials over  $k$ equipped with the  symbolic product 
$\star$  is isomorphic to the $\F_q$-algebra  $k\{\tau_q\}$, 
the isomorphism being the following:
$$P(T)=\sum_{j=0}^{d}\alpha_jT^{q^j}\mapsto  P(\tau_q)=\sum_{j=0}^{d}\alpha_j\tau_q^j.$$
We say that $P(\tau_q)$ is monic (resp. separable) if 
$P(T)$ is monic (resp. separable). Notice then that $P(\tau_q)=\sum_{j=0}^{d}\alpha_i\tau_q^{j}$ is separable if and only if $\alpha_0\neq 0$. Moreover
$$\Deg P(T)=q^{\Deg P(\tau_q)}.$$

\noindent
The set of roots of a $q$-polynomial over $k=\F_r$ has special properties.
\begin{Lemme}
\label{subspace}
Assume that $\F_q\subset \F_r\subset \F_s.$
\begin{enumerate}
\item Let $P(T)\in\F_r[T]$ be a non-zero $q$-polynomial over $\F_{r}$ 
and let $\F_{s}$ be a finite extension 
of $\F_{r}$ containing all the roots of $P$. Then each root of 
$P(T)$ has the same multiplicity, which is either $1$ or a power of $q$, 
and the set of roots of $P(T)$ is a $\F_q$-subspace of $\F_{s}$.
\item 
Let $H$ be a $\F_q$-subspace of $\F_s$. Then
$P_H(T)=\prod_{a\in H}(T-a)$
is a $q$-polynomial over $\F_s$. Moreover, if $H$ is stable by $\Gal(\F_s/\F_r)$, 
then $P_H(T)$ is a $q$-polynomial over $\F_r$.
 \end{enumerate}
\end{Lemme}
\begin{proof} see \cite[Theorems 3.50 and  3.52. p. 103]{Lidl}. The last assertion is clear.
\end{proof}
If $H$ is a $\F_q$-subspace spanned by $\{w_1,\ldots,w_n\}\subset k$, we set
$H=\langle w_1,\ldots,w_n\rangle_q$.

\begin{Defi}
Let  $\F_q\subseteq k$ and $\{w_1,\ldots,w_d\}\subset k$. We define the {\em Moore determinant of $\{w_1,\ldots,w_d\}$ over $\F_q$} to be
$$\Delta(w_1,\ldots,w_d)=\left\vert
\begin{array}{ccc}
w_{1}&\cdots&w_d\\
w_{1}^q&\cdots&w_d^q\\
\vdots&&\vdots\\
w_{1}^{q^{d-1}}&\cdots&w_{d}^{q^{d-1}}\\
\end{array}
\right\vert\,.
$$
\end{Defi}
We refer to  \cite[page 8]{Goss} for the computation of a Moore determinant.

\begin{Lemme}
\label{moore}
 Let $k$ be a field such that $\F_q\subseteq k$ and
let $H$ be a $\F_q$-subspace of $k$ of dimension $d$. Then $\{w_1,\ldots,w_d\}\subset H$ is a basis of $H$ over $\F_q$  if and only if $\Delta(w_1,\ldots,w_d)\neq 0$.
\end{Lemme}
\begin{proof}
 \cite{Goss}, Corollary 1.3.4.
\end{proof}

\subsection{Equations of subextensions of an Artin-Schreier type extension}
The following result is an application of  \cite{Goss}, Proposition 1.3.5.

\begin{Propo}
\label{descente}
Set $q=p^n$, $r=q^e$ and let $A(T)\in \F_r[T]$ be a monic non-zero  separable $q$-polynomial over $\F_r$ of degree $q^d$,  $1\leq d\leq e$, having all its roots in an extension $\F_s$   of $\F_r$. 
The set of roots of $A(T)$, denoted by $\cal A$, 
is a $\F_q$-vector space of dimension $d$.
Assume that  ${\cal A}$ and all its $\F_q$-subspaces are stable under $\Gal(\F_s/\F_r)$.  
\begin{enumerate}
\item
Let $(w_1,\ldots,w_d)$ be a basis of $\cal A$ over $\F_q$. For all $i$, $1\leq i<d$, 
set
$$ A_i (T)=\prod_{a\in H_i }(T-a),\,\mbox{ where } H_i=\langle w_1,\ldots,w_i\rangle_q.$$
Then $A_i(T) $ is a monic separable $q$-polynomial over $\F_r$ of degree $q^i $ 
and there exists a unique monic separable $q$-polynomial over $\F_r$ of degree $q^{d-i }$, 
$M_i(T) $, such that $A(T)=(M_i \star A_i)(T) $. More precisely
$$M_i (T)=\prod_{a\in \bar H_i }(T-a),\, \mbox{ where } 
\bar H_i=\langle A_i (w_{i +1}),\ldots, A_i (w_{d}) \rangle_q.
$$
\item Let $E/\F_{r}$ be an algebraic function field
and set $G=E(z)$, where $A(z)=u\in E$.
Assume moreover that there exists a place $\wp$ of $E$ such that $\nu_\wp(u)=-m$, with $m>0$ and $\gcd(m,p)=1$. 
For all $i $, $1\leq i <d$, we consider the subfield $G_i$ of $G$ defined by 
$$G_{i} =E(t_{i} ), \mbox{ with }t_{i} =A_{d-i} (z) \mbox{ and }M_{d-i} (t_{i})-u=0.$$
Then the full constant field of $G$ and $G_{i} $ is $\F_r$ and $[G_{i} :E]=q^{i }$.
\end{enumerate}
\end{Propo}

\begin{proof}
\begin{enumerate}
 \item By Lemma \ref{subspace}, $\cal A$ is a $\F_q$-vector space.
 The polynomial $A_i(T)=\prod_{a\in H_i}(T-a)$ is a monic factor of $A(T)$ and its degree is equal to $q^i$. By Lemma \ref{subspace} and
since $H_i$ is a $\F_q$-vector space,  $A_i(T)$ is   a $q$-polynomial over $\F_s$. But, since $A_i(T)$ is stable under $\Gal(\F_{s}/\F_r)$, it is  a 
$q$-polynomial over $\F_r$.
Set $\sigma_j=A_i(w_j)$, for all $j=1,\ldots, d$ and write $A_i(T)=\sum_{m=0}^{i}b_mT^{q^m}$, with $b_i=1$ and $b_m\in\F_r$ for all $1\leq m\leq i$.
Since $\sigma_j=0$, for all $j=1,\ldots, i$, we have
$\Ker (A_i)=  H_i$ and $\im(A_i)=<\sigma_{i+1},\ldots,\sigma_{d}>_q$, where we denote by $A_i$ the $\F_q$-linear map from ${\cal A}$ into $\F_{s}$ associated to $A_i(T)$.
In particular, $\sigma_{i+1},\ldots,\sigma_{d}$ are $t=d-i$  vectors of $\F_{s}$  which are independent over $\F_q$.
We look for a polynomial
$M_i(T)=\sum_{m=0}^{d-i}c_mT^{q^m}$, with $c_{d-i}=1$ and $c_m\in \F_r$, such that 
$A(T)=M_i(A_i(T))$. But this  is equivalent to $M_i(\sigma_j)=0$, for all $j=i+1,\ldots, d$.
Thus the $c_m$'s are solutions of the following system of $t$ equations and $t$ indeterminates:
$$
\left\{
\begin{array}{ccccccccc}
\sigma_{i+1}^{q^{t-1}}X_{t-1}&+&\sigma_{i+1}^{q^{t-2}}X_{t-2}&+
&\cdots&+&\sigma_{i+1}X_0&=&-\sigma_{i+1}^{q^{t}}\\
\vdots&&&&&\\
\sigma_{d}^{q^{t-1}}X_{t-1}&+&\sigma_{d}^{q^{t-2}}X_{t-2}&+
&\cdots&+&\sigma_{d}X_0&=&-\sigma_{d}^{q^{t}}.\\
\end{array}
\right.
$$
By Lemma \ref{moore}, the Moore determinant over $\F_q$ of this system is non zero if and only if $\sigma_{i+1},\ldots,\sigma_{d}$ are linearly independent over $\F_q$, which is the case.
Thus the preceding system has a unique solution in $\F_{s}^t$. Since $\cal A$ and $H_i$ are stable under $\Gal(\F_{s}/\F_r)$, we obtain that $M_i(T)\in \F_r[T]$. Notice that $A_i(T)=\sum_{m=0}^{i}b_mT^{q^m}$ and $M_i(T)=\sum_{m=0}^{d-i}c_mT^{q^m}$ are separable polynomials, since the product $c_0b_0$ is equal to the coefficient of $T$ in $A(T)$ which is non-zero. 

\item Consider the constant field extension $L/\F_s$ of the function field $E/\F_r$
and set $F=G\otimes_{\F_r}\F_s$. Then $F=L(z)$, with $A(z)-u=0$.  First, we show that the full constant field of $F$ is $\F_s$.
 Let $\mathfrak q$ be a place of $L/\F_s$ above $\wp$. Since the extension $L/\F_s$ of $E/\F_r$ is unramified, we have $\nu_{\mathfrak q}(u)=\nu_\wp(u)=-m$, with $\gcd(m,p)=1$.
Thus the polynomial $A(T)-u\in E[T]\subset L[T]$ is absolutely irreducible. We also deduce from this fact that the full constant field of $G$ is $\F_r$.
Moreover the  additive separable polynomial $A(T)$ has all its roots in $\F_s$ and  Proposition III.7.10 of \cite{stich} says that the extension $F/L$ 
is Galois of degree $q^d$ and
$\Gal(F/L)\simeq{\cal A}$. We identify $\Gal(F/L)$ with ${\cal A}$.

To the sequence of subgroups $H_1\subset H_2\subset\ldots \subset H_{d-1}\subset {\cal A}$, 
we associate the sequence of subfields $L_i/\F_s$ of $F/\F_s$ 
$$L\subset L_1\subset \ldots \subset L_{d-1}\subset F,$$
where $L_i=F^{H_{d-i}}$ is the fixed field of $H_{d-i}\subseteq\Gal(F/L)$. By Galois theory, 
$[F:L_i]=\ord H_{d-i}=q^{d-i}$ and thus $[L_i:L]=q^i$. 
The element $t_i=A_{d-i} (z)=\prod_{a\in H_{d-i}}(z-a)\in \F_r[z]$ is in $F$ 
and it is fixed by $H_{d-i}$, then $L(t_i)\subseteq  L_i$. Moreover, 
the minimal polynomial of $z$ over $L(t_i)$ divides $A_{d-i}(T)-t_i$ 
which is of degree $q^{d-i}$. 
Thus $L_i=L(t_i)$ and the minimal polynomial of $z$ over $L_i$ is $A_{d-i}(T)-t_i$.
Finally,  $[L_i:L]=q^i$ and the minimal polynomial of $t_i$ over $L$ is $M_{d-i}(T)-u$.
Since $t_i\in \F_r[z]$ is an element of $G$, the field $G_i=E(t_i)$ is a subfield of $G$, 
and the minimal polynomial of $t_i$ over $E$ is $M_{d-i}(T)-u$ again, 
since $M_{d-i}(T)\in \F_r[T]$ and $u\in E$.
Clearly, $L_i=G_i\otimes_{\F_r}\F_s$, which proves that $L_i/\F_s$ 
is the constant field extension of $G_i/\F_r$.
To summarize, we have
$$\xymatrix{
G=E(z) \ar@{-}[r] & F=L(z) \\
G_i=E(t_i) \ar@{-}[u]^{q^{d-i}} \ar@{-}[r] & L_i =L(t_i)\ar@{-}[u]_{q^{d-i}}\\
E \ar@{-}[u]^{q^{i}} \ar@{-}[r] & L\ar@{-}[u]_{q^{i}} \\
\F_r \ar@{-}[u] \ar@{-}[r] &\F_s\ar@{-}[u]
}$$
\end{enumerate}
\end{proof}

We apply the previous result to the polynomial $A(T)=T^q+T$. 

\begin{Coro}
\label{descenteAS}
 Set $q=p^n$ and $P(T)=T^q+T$.
 \begin{enumerate}
 \item[\rm (i)]  The set 
 ${\cal P}\subset \F_{q^2}$  of roots of the separable $p$-polynomial 
 $P(T)$ is a $\F_p$-vector space of dimension $n$, 
 which is stable under $\Gal(\F_{q^2}/\F_q)$. 
 Let $(w_1,\ldots,w_n)$ be a basis of $\cal P$ over $\F_p$.
 \item[\rm (ii)] For all $i$, $1\leq i<n$, we set
 \begin{equation}
 \label{Pi}
 P_i (T)=\prod_{a\in H_i }(T-a),\,\mbox{ where } H_i=\langle w_1,\ldots,w_i\rangle_p.
 \end{equation}
Then $P_i $ is a monic separable $p$-polynomial over $\F_q$ of degree $p^i $ 
and there exists a unique monic separable $p$-polynomial over $\F_q$ of degree $p^{n-i }$, 
$M_i $, such that $P(T)=(M_i \star P_i)(T) $. More precisely
 \begin{equation}
 \label{Mi}
M_i (T)=\prod_{a\in \bar H_i }(T-a),\, \mbox{ where } \bar H_i=\langle P_i (w_{i +1}),\ldots, P_i (w_{n} \rangle_p.
\end{equation}
\item[\rm (iii)] Let $L/\F_{q^2}$ be an algebraic function field and consider an Artin-Schreier type extension 
$F/\F_{q^2}$ of $L/\F_{q^2}$ defined by $F=L(z)$ with $P(z)=u\in L$.
We assume that 
there exists a place $\wp$ of $L$ such that $\nu_\wp(u)=-m$, with $m>0$ and $\gcd(m,p)=1$. 
Assume moreover that $L/\F_{q^2}$ can be reduced over $\F_q$ and let $E/\F_q$ be 
the reduced field. We assume that $u\in E$ and set $G=E(z)$, where $P(z)=u$.
For all $i $, $1\leq i <n$, we consider the sub-function field $G_i$ of $G$ defined by 
$$G_{i} =E(t_{i} ), \mbox{ with }t_{i} =P_{n-i} (z) \mbox{ and }M_{n-i} (t_{i})-u=0.$$
Then the full constant field of $G$ and $G_{i} $ is $\F_q$, $[G_{i} :E]=p^{i }$ and
$E\subset G_{1}\subset \ldots \subset G_{n-1}\subset G$.
\end{enumerate}
\end{Coro}

\begin{proof}
\begin{enumerate}
 \item[\rm (i)]  See Example \ref{$Tq+T$}.  
We have seen that 
${\cal P}$ and also any $\F_p$-subspace of $\cal P$ is stable under $\Gal(\F_{q^2}/\F_q)$.
Let $(w_1,\ldots,w_n)$ be a basis of ${\cal P}$ and let $H_i$ be the $\F_p$-subspace spanned by $(w_1,\ldots,w_i)$, for $1\leq i<n$.

 \item[\rm (ii)] Apply assertion 1 of Proposition \ref{descente}.

 \item[\rm (iii)] We have assumed that $u\in E$, where $E/\F_q$ is 
 the constant restriction (or descent) of  
 $L/\F_{q^2}$ and that there exists a place $\wp$ of $L$ 
 such that $\nu_\wp(u)=-m$, with $m>0$ and $\gcd(m,p)=1$. 
 Thus there exists a unique place $\mathfrak p$ of $E/\F_q$ which lies under $\wp$ and such that $\nu_{\mathfrak p}(u)=-m$.
Then we apply assertion 2 of Proposition \ref{descente}.
Assume $n\geq 2$ and let us show that $G_i\subset G_{i+1}$ for all $1\leq i<n-1$.
We have $[G:G_i]=p^{n-i}$ and $P_{n-i}(z)=t_i$ with $\Deg P_{n-i}(T)=p^{n-i}$ for all $1\leq i<n$.
The equation of $G$ over $G_i$ is then $P_{n-i}(z)=t_i$.
We can apply assertion 2 of Proposition \ref{descente} recursively to obtain $G_1$ such that 
$E\subset G_1\subset G$, then $G_2$ such that  $G_1\subset G_2\subset G$ and so on. 
But  also, 
applying assertion 1 in \cite{Goss}, Proposition 1.3.5, we have for all $1< j<n$
$$
P_{j}(T)=P_{j-1}(T)^p-W_jP_{j-1}(T), 
$$
where $W_j=P_{j-1}(w_j)^{p-1}=\prod_{a\in H_{j-1}}(w_j-a)^{p-1}$.
Moreover $W_j$ is in $\F_q$, since
$$W_j^q=\prod_{a\in H_{j-1}}(w_j^q-a^q )^{p-1}
=\prod_{a\in H_{j-1}}(-w_j+a)^{p-1}=W_j.$$
Thus $G_{i+1}=G_{i}(t_{i+1})$, where
\begin{equation}
\label{ti}
t_{i+1}^p-W_{n-i}t_{i+1}=t_{i}.
\end{equation}
\end{enumerate}
\end{proof}

\noindent
We can apply the preceding result to the Hermitian function field $F/\F_{q^2}$. The fact that the Hermitian function field has a descent from $\F_{q^2}$ to $\F_q$ is obvious since its equation is defined over $\F_q$. Moreover, 
if $p=2$, the reduced field of $F$ over $\F_q$ is a Galois extension of the rational function field $\F_q(x)$, thus the existence of the intermediate steps is also obvious (see \cite{ballro}).

\begin{Coro}
\label{descenteGS} 
Set $G_1=\F_q(x)$ and
consider the descent over $\F_q$ of the Hermitian function field, denoted by
$G_2/\F_q=\F_q(z,x)/\F_q$, with $z^q+z=x^{q+1}$.
For all $i $, $1\leq i <n$, we consider the subfield $G_{1,i}$ of $G_2$ defined by 
$$G_{1,i} =\F_q(x,t_{i} ), \mbox{ with }t_{i} =P_{n-i} (z)\mbox{ and }M_{n-i} (t_{i})-x^{q+1}=0,$$
where $P_i$ and $M_i$ are defined by $(\ref{Pi})$ and $(\ref{Mi})$ respectively.
Then the full constant field of $G_{1,i} $ is $\F_q$, $[G_{1,i} :\F_q(x)]=p^{i }$, 
the genus of $G_{1,i} /\F_q$ is $g_i =\frac{q(p^i -1)}2$ and $G_1\subset G_{1,1}\subset \ldots \subset G_{1,n-1}\subset G_2$.
\end{Coro}

\begin{proof} Since the equation of the Hermitian function field $F_2/\F_{q^2}$ 
is defined over $\F_q$, it has a descent over $\F_q$, say $G_2/\F_q$.
Consider the rational function fields  $F_1/\F_{q^2}=\F_{q^2}(x)/\F_{q^2}$
and $G_1/\F_{q}=\F_{q}(x)/\F_{q}$. The infinite place $\wp$ of $F_1/\F_{q^2}$ 
is fully ramified in 
$F_2/\F_{q^2}$ and $\nu_\wp(x^{q+1})=-(q+1)$.
Let us denote by $\mathfrak p$ the infinite place of $G_1/\F_{q}$. 
Then $\mathfrak p$ is fully ramified in $G_2/\F_{q}$.
Since $G_{1,i}$ is a sub-extension of $G_2$, $\mathfrak p$ is 
also fully ramified in $G_{1,i}$ and the result follows. 
Notice that, if $p=2$, the extension $G_2/G_{1,i}$ is Galois, otherwise it is not Galois.
To compute the genus of $G_{1,i}/\F_q$, we apply \cite[Prop.  III.7.10]{stich},
considering the constant field extension $F_{1,i}/\F_{q^2}$ of $G_{1,i}/\F_q$. 
The genus of $G_{1,i}/\F_q$ and $F_{1,i}/\F_q$ are equal and $F_{1,i}=\F_{q^2}(x,t_{i})$, 
with 
$M_{n-i}(t_{i})-x^{q+1}=0$. Since $M_{n-i}(T)$ is an additive 
polynomial of degree $p^i$ which has all its roots in 
$\F_{q^2}$ and since $\nu_{\wp}(x^{q+1})=-(q+1)$,
Proposition  III.7.10 of \cite{stich} shows that 
the genus of $F_{1,i}/\F_{q^2}$ is  equal to 
$g_i=\frac{q(p^i-1)}2$.
\end{proof}
Notice that the constant field extension $F_{1,i}/\F_{q^2}$ of each $G_{1,i}/\F_q$ 
is a subfield of the Hermitian function field and thus is maximal.
The number of rational places of each $G_{1,i}/\F_q$ is then equal to $q+1$.

\begin{ex} 
Assume $p=2$ and $q=p^n$. 
Then ${\cal P}=\F_q$ and a $\F_p$-basis of $\cal P$ 
is just a basis of $\F_q$ over $\F_p$. If $w$ is a generator of $\F_q$ over $\F_2$, 
we consider the groups $H_{i+1}=\langle 1,w,\ldots,w^i\rangle_2$, for $0\leq i<n-2$.
We have seen previously that the $p$-polynomials $Q(T)=T^2+T$ and  
$M(T)=T^{2^{n-1}}+\cdots+T$ are symbolic factors of 
$P(T)$ and $P(T)=T^q+T=M(Q(T))=Q(M(T))$. Moreover, 
$M(T)$ defines the trace map from $\F_{2^n}$ to $\F_2$.
The element
$t_{n-1}=Q(z)=z(z+1)$ corresponds to $H_1$ and
this gives step $n-1$ 
$$G_{1,n-1}=G_1(t_{n-1}),\, \mbox{ where  }
M(t_{n-1})+x^{q+1}=0.$$
If $n=2$, we are done.
Assume $n>2$ is even. Then
$H_{n-1}=\{a\in\F_q,\,M(a)=0\}$ is a subspace of $\cal P$ of order $2^{n-1}$
which contains $H_{1}$.
This gives step $1$, $t_1=\prod_{a\in H_{n-1}}(z+a)$, $G_{1,1}=G_1(t_{1})$ and the minimal polynomial of $t_1$ over $G_1$ is $Q(T)+x^{q+1}=0.$ 
Notice that, if $n$ is odd, we cannot obtain  step $1$ by this method.

\begin{enumerate}
\item $q=2^2$. 
$t_1=z(z+1)$, and $G_{1,1}=\F_q(x,t_1)$, with 
$M_1(t_1)=x^{q+1}$, where $M_1(T)=T^2+T$.

\item $q=2^3$. Let $w$ be a generator of $\F_q$,  with $w^3+w+1=0$.
The two intermediate function fields $G_{1,i}/\F_q$ are $G_{1,i}=\F_q(x,t_i)$, $i=1,\,2$, where
\begin{enumerate}
\item $t_2=z(z+1)$,  $M_1(T)=T^4+T^2+T$,
\item $t_1=t_1(z+w)(z+w+1)$, $M_2(T)=T^2+w^3T$.
\end{enumerate}

\item $q=2^4$.  Let $w$ be a generator of $\F_q$,  with $w^4+w+1=0$.
\begin{enumerate}
\item $t_3=z(z+1)$,  
$M_1(T)=T^8+T^4+T^2+T$,
\item $t_2=t_3(z+w)(z+w+1)$, $M_2(T)=T^4+w^{10}T^2+w^{10} T$,
\item $t_1=M_1(z)$, $M_3(T)=T^2+T$,
\end{enumerate}

\item $q=2^5$.  Let $w$ be a generator of $\F_q$,  with $w^5+w^2+1=0$.
\begin{enumerate}
\item $t_4=z(z+1)$,  
$M_1(T)=T^{16}+T^8+T^4+T^2+T$,
\item $t_3=t_4(z+w)(z+w+1)$, $M_2(T)=T^8 + w^{26}T^4 + w^{16}T^2 + w^{12}T$,
\item $t_2=t_3\prod_{a\in\{w^2,w^{19}\}}(z+a)(z+a+1)$,
$M_3(T)=T^4 + w^{29}T^2 + w^6T$,
\item 
$t_1=t_2\prod_{a\in\{w^3,w^{6},w^{12},w^{20}\}}(z+a)(z+a+1)$, $M_4(T)=T^2 + w^{2}T$.
\end{enumerate}
\end{enumerate}
Let us apply formula $(\ref{ti})$ in the last case $q=2^5$:
$$\xymatrix{
z^{32}+z=x^{33}&\ar@{-}[ddddd]_{2^5}G_2=G_1(z)&G_2=G_{1,4}(z) &z^2+z=t_4 \\
&&G_{1,4}=G_{1,3}(t_4) \ar@{-}[u]_2& t_4^2+w^{19}t_4=t_3\\
&&G_{1,3}=G_{1,2}(t_3) \ar@{-}[u]_2& t_3^2+w^{6}t_3=t_2\\
&&G_{1,2}=G_{1,1}(t_2) \ar@{-}[u]_2& t_2^2+w^{4}t_2=t_1\\
&&G_{1,1}=G_{1}(t_1) \ar@{-}[u]_2&  t_1^2+w^{2}t_1=x^{33}\\
&G_1&G_1=\F_{32}(x)\ar@{-}[u]_2&\\
}$$

Assume  that $p\geq 3$. Then, setting
${\cal N}_q({\cal P}):=\norm_{\F_{q^2}/\F_q}({\cal P}\setminus\{0\}),$
the factorization of $P$ in $\F_q[T]$ is
$$P(T)=T\prod_{a\in {\cal N}_q({\cal P})}(T^2+a)\,.$$
For instance, for $p=3$ we have:
\begin{enumerate}
\item $q=3^2$.  Let $w$ be a generator of $\F_q$,  with $w^2+2w+2=0$. Then
$$T^q+T=T\prod_{a\in\{w,w^3,w^5,w^7\}}(T^2+a).$$
Set $t_1=z(z^2+w)$, then the minimal polynomial of $t_1$ over $G_1=\F_q(x)$ is
$$M_1(T)-x^{q+1} \hbox{, where } M_1(T)=T^3+w^7T\,.$$
\item $q=3^3$.  Let $w$ be a generator of $\F_q$,  with $w^3+2w+1=0$. Then
$T^q+T=T\prod_{i=0}^{12}(T^2+w^{2i}).$
\begin{enumerate}
\item   $t_2=z(z^2+1)$, $M_1(T)=T^9+2T^3+T\,.$
\item  $t_1=t_2(z^2+w^2)(z^2+w^6)(z^2+w^{18})$, $M_2(T)=T(T^2+1)$.
\end{enumerate}
\end{enumerate}
\end{ex}

\section{Applications}\label{trois}

In this section $q=p^n$ is an arbitrary prime power, $K=\F_{q^2}$ and $k=\F_q$. 

\subsection{Descent of a completed Garcia-Stichtenoth's tower}
Let us recall the definition of a separable tower of function fields 
over a field $K$ (cf. \cite{stich}).

\begin{Defi}
A separable tower ${\cal T}$ over $K$ 
is a sequence $(F_i)_{i\geq 1}$ of algebraic function fields of one variable over $K$
such that
$${\cal T}:\,F_1 \subsetneq F_2 \subsetneq \cdots \subsetneq F_i \subsetneq \cdots ,$$
and that each extension $F_i/F_{i-1}$ is finite separable.
\end{Defi}

 \medskip

We consider the Garcia-Stichtenoth's tower ${\cal T}_{GS}$ over $\F_{q^2}$ 
constructed in \cite{gast}. Recall that this tower
is defined recursively in the following way.
We set $F_1=\F_{q^2}(x_1)$ the rational function field over  $\F_{q^2}$,
and for $i \geq 1$ we define 
$$F_{i+1}=F_i(z_{i+1}),$$
where $z_{i+1}$ satisfies the equation
$$z_{i+1}^q+z_{i+1}=x_i^{q+1},$$
with
$$x_i=\frac{z_i}{x_{i-1}} \hbox{ for } i\geq 2.$$

We consider the completed Garcia-Stichtenoth's tower ${\cal T}_c$ over $\F_{q^2}$ 
studied in \cite{ball3} 
obtained from ${\cal T}_{GS}$
by adjonction of intermediate steps. Namely
we have
$${\cal T}_c:\, F_{1,0}\subset \cdots  
\subset F_{i,0} \subset F_{i,1}\subset\cdots \subset F_{i,s} 
\subset \cdots \subset F_{i,n-1}
\subset F_{i+1,0}\subset \cdots $$
where the steps $F_{i,0}$ are the steps $F_i$ of the Garcia-Stichtenoth's tower
and where $F_{i,s}$ ($1 \leq s \leq n-1$) are the intermediate steps.
 
Each extension $F_{i,s}/F_i$ is Galois of degree $p^s$,
with full constant field $\F_{q^2}$. More precisely, 
each extension $F_{i,s}/F_i$ is an elementary abelian $p$-extension.
Namely, the Galois group  $\Gal(F_{i,s}/F_i)$ is isomorphic to
$\left ( \Z/p\Z \right)^s $.


Let us mention the following well known result (cf. \cite{stich}, Proposition III.7.10):

\begin{Propo}\label{util}
Consider an algebraic function field $F/K$, a linearized separable polynomial $A(T)\in K[T]$ of degree $p^s$ which has all its roots in $K$ and $u \in F$. Assume that the polynomial $A(T)-u$  is absolutely irreducible and set $E=F(z)$ with $A(z)=u$.
Then $E/F$ be an elementary abelian $p$-extension of degree $p^s$ and
 the Galois group of $E/F$ is 
$$\Gal(E/F)=\{z\mapsto z+\alpha~|~A(\alpha)=0 \hbox{ and } \alpha\in K  \}.$$
\end{Propo}

\begin{Theo}
\label{tourc} Set $q=p^n$.
The descent of the definition field 
of the tower ${\cal T}_c$  from $\F_{q^2}$ to $\F_q$ is possible.
More precisely, there exists a tower ${\cal R}_c$ over $\F_q$
given by a sequence 
$${\cal R}_c:\,G_{1,0}\subset \cdots  
\subset G_{i,0}\subset G_{i,1}\subset\cdots \subset G_{i,s} 
\subset \cdots \subset G_{i,n-1}
\subset G_{i+1,0}\subset \cdots $$
of algebraic function fields with full constant field $\F_q$
such that, for all $i\geq 1$ and $1\leq s<n$ 
$$F_{i,s}=G_{i,s} \otimes_{\F_q} \F_{q^2}.$$
\end{Theo}

\begin{proof}
For any
integer $i\geq 1$, let $E=F_{i,s}$ with $1 \leq s \leq n-1$
be an intermediate step between $F_{i,0}$ and $F_{i+1,0}$.
It is sufficient to prove that the Galois group
$\Gamma=\Gal(\F_{q^2}/\F_q)$ acts by inner automorphisms on
the Galois group $H=\Gal(F_{i+1,0}/E)$. We know that $F_{i+1,0}/F_{i,0}$
is an elementary abelian $p$-extension of Artin-Schreier type of degree $q=p^n$.
Hence, by Proposition \ref{util} the Galois group $G$ of the extension $F_{i+1,0}$
over $F_{i,0}$ is the set: 
$$\Gal(F_{i+1,0}/F_{i,0})= 
\{z_{i+1}\mapsto z_{i+1}+\alpha ~|~ \alpha^q+\alpha=0,~ \alpha \in \F_{q^2}\},$$ 
which is isomorphic to the additive group $(\Z/p\Z)^n$.
The group $\Gamma$ is constituted by the identity and the Frobenius
automorphism $\tau_q$. The automorphism $\tau_q$ acts on
the Galois group $G$ by the transformation of $\alpha$ into $-\alpha$.
Hence, the subgroup $H$ of $G$ is let invariant by $\Gamma$.
The result follows from Theorem \ref{princ} and Remark \ref{rprinc}. 
\end{proof}

We can prove the preceding result in an effective way.
In the following theorem, we use the notations $x_i$ and $z_i$ introduced in the definition of  the tower ${\cal T}_{GS}$.
\begin{Theo} Set $q=p^n$. The descent of the definition field 
of the tower ${\cal T}_c$  from $\F_{q^2}$ to $\F_q$ is possible.
The reduced tower over $\F_q$ is ${\cal R}_c=(G_{i,s},\, i\geq 1,\, 0\leq s<n)$, where $G_{1,0}=\F_q(x_1)$ and, for all $i\geq 1$,
the explicit equations of the function fields $G_{i,s}/\F_q$ are the following:
$$\begin{array}{ll}
G_{i+1,0}=G_{i,0}(z_{i+1}),\,&z_{i+1}^q+z_{i+1}=x_i^{q+1}\\
G_{i,s}=G_{i,0}(t_{i,s}),\,&t_{i,s}=P_{n-i}(z_{i+1}),\, M_{n-i}(t_{i,s})=x_i^{q+1},\\
\end{array}$$
where $P_j(T)$ and $M_j(T)$, for $1\leq j<n$, are the $p$-polynomials over $\F_q$
defined 
by $(\ref{Pi})$ and $(\ref{Mi})$.
\end{Theo}
\begin{proof} This follows readily from  Corollary \ref{descenteAS}, since each extension $F_{i+1,0}/F_{i,0}$ in tower ${\cal T}_c$ is of Artin-Schreier type and the infinite place of $F_{1,0}/\F_{q^2}$ is fully ramified in $F_{i,0}/\F_{q^2}$, for all $i\geq 2$. Notice that the intermediate steps in the first stage $G_{2,0}/G_{1,0}$ are given in Corollary  \ref{descenteGS}.
\end{proof}

\subsection{On the bilinear complexity of the multiplication}

We denote by $m$ the ordinary multiplication  in the finite field $\F_{q^n}$
of characteristic $p$. 
This field will be considered as a $\F_q$-vector space. 
The multiplication $m$ is a bilinear map from $\F_{q^n} \times \F_{q^n}$ into $\F_{q^n}$, 
thus it corresponds to a linear map $M$ from the tensor product 
$\F_{q^n} \bigotimes \F_{q^n}$ over $\F_q$ into $\F_{q^n}$. 
One can also represent 
$M$ by a tensor $t_M \in \F_{q^n}^*\bigotimes \F_{q^n}^* \bigotimes \F_{q^n}$ 
where $\F_{q^n}^*$ denotes the dual of $\F_{q^n}$ over $\F_q$. 
Hence the product of two elements $x$ and $y$ of $\F_{q^n}$ 
is the convolution of this tensor with $x \otimes y \in \F_{q^n} \bigotimes \F_{q^n}$. 
The tensor rank $\mu_{q}(n)$ of $t_M$ is called the bilinear complexity of 
multiplication in $\F_{q^n}$ over $\F_q$. It corresponds 
to the minimum possible number of summands in any tensor decomposition.  

\medskip

The following theorem, proved in \cite{ballro},
provides an estimation of the bilinear complexity $\mu_{q}(n)$
under the assumption that there exists well-fitted  
function fields.

\begin{Theo}\label{Chud1}
Let $q$ be a prime power and $n>1$ an integer. If there exists a function
field $F/\F_q$ of genus $g(F)$ with $N_1(F)$ places of degree $1$
and $N_2(F)$ places of degree $2$ such that

\medskip

\quad 1) there exists a non-special divisor of degree $g(F)-1$ ,

\par \quad 2) $2g(F)+1 \leq q^{\frac{n-1}{2}}(q^{\frac{1}{2}}-1)$,

\par \quad 3) $N_1(F)+2N_2(F) > 2n+2g(F)-2$,

then

$$\mu_q(n) \leq 3n+3g(F).$$
\end{Theo}
 
The existence of a descent for the tower ${\cal T}_c$
and the results of \cite{ball4}
on the existence of a non-special divisor a degree $g-1$
for all steps of the tower ${\cal T}_c$ when $q \geq 4$, make it legitimate to apply
Theorem \ref{Chud1}.
We can now use for any characteristic $p$ 
the method developped in \cite{ballro} for the
particular case $p=2$.  We obtain :

\begin{Theo}\label{comp}
For $q \geq 4$ and for any integer $n$ 
we have
$$
\mu_q(n) \leq 
  3(1+\frac{p}{q-3})n,
$$ 
and
$${\cal M}_q=\limsup_{n\rightarrow +\infty} 
\frac{\mu_q(n)}{n}\leq 3(1+\frac{p}{q-3}).$$
\end{Theo}

Let us remark that in the case $q=p \geq 5$
the bound given in \cite{ballchau}: 
$$
\mu_p(n) \leq 
  3(1+\frac{4}{p-3})n
$$
is better. 
The asymptotic bound given in Theorem \ref{comp}
is better than the one  given in \cite{shtsvl} and \cite{ballro}.
However for $q=p \geq 5$, the best asymptotic bound is the one given
in \cite{ballchau}.

\medskip

 {\bf Acknowledgements:} the authors wish to thank J.M. Couveignes 
 for many valuable discussions.


\bigskip

{\bf S. Ballet} , 
Laboratoire de G\'eom\'etrie Alg\'ebrique et Applications \`a la
Th\'eorie de l'Information,
Universit\'e de la Polyn\'esie Fran\c caise, BP 6570, 98702 Faa'a,
Tahiti, Polyn\'esie Fran\c caise (France).
\par 
e-mail: ballet@upf.pf.

\medskip

{\bf D. Le Brigand} , Institut de Math\'ematiques de Jussieu,
Universit\'e Pierre et Marie Curie, Paris VI, Campus Chevaleret,
175 rue du Chevaleret, F75013 Paris.
\par
e-mail: Dominique.LeBrigand@math.jussieu.fr.

\medskip

{\bf R. Rolland} , C.N.R.S. Institut de Math\'{e}matiques
de Luminy, case 930, F13288 Marseille cedex 9 (France).
\par
e-mail: rolland@iml.univ-mrs.fr.


\begin{thebibliography}{0}

\bib{ball2}{S. Ballet,}
{Curves with Many Points and Multiplication Complexity in Any Extension of $\F_q$.}
{Finite Fields and Their Applications, 5 (1999), 364-377.}

\bib{ball3}{S. Ballet,}
{Low Increasing Tower of Algebraic Function
Fields and Bilinear Complexity
of Multiplication in Any Extension of $\F_q$.}
{Finite Fields and Their Applications, 9 (2003), 472-478.}

\bib{ball4}{S. Ballet, D. Le Brigand,}
{On the existence of Non-Special Divisors
of degree $g$ and $g-1$ in Algebraic Function Fields over $\F_q$.}
{Preprint.}

\bib{ballchau}{S. Ballet, J. Chaumine,}
{On the bounds of the bilinear complexity of multiplication in some
finite fields.}
{To appear in AAECC.}

\bib{ballro}{S. Ballet, R. Rolland,}
{Multiplication Algorithm in a Finite Field and Tensor Rank of the Multiplication.}
{Journal of Algebra, 272 (2004), 173-185.}

\bib{Deo1}{V. Deolalikar,}
{Determining irreducibility and ramification groups for 
an additive extension of the rational function field.}
{Journal of Number Theory, 97 (2002), no. 2, 269-286.}


\bib{gast2}{A. Garcia, H. Stichtenoth,}
{Elementary Abelian $p$-Extensions of Algebraic Function Fields.}
{Manuscripta Mathematica, 72 (1991), 67-79.}

\bib{gast}{A. Garcia, H. Stichtenoth,}
{A tower of Artin-Schreier extensions of function 
fields attaining the Drinfeld-Vladut bound.}
{Inventiones Mathematicae, 121 (1995), 211-222.}

\bib{stiGar}{A. Garcia, H. Stichtenoth,}
{A class of polynomials over finite fields.}
{Finite Fields and Their Applications, 5 (1999), no. 4, 424-435.}

\bib{Goss}{D. Goss,}
{Basic structures of function fields arithmetic.}
{Springer 1991.}


\bib{Lidl}{R. Lidl, H. Niederreiter,}
{Introduction to finite fields and their applications.}
{Cambridge Univ. Press, 1994 (revised version)}.

\bib{serre}{J.-P. Serre,}
{Groupes Alg\'ebriques et Corps de Classes.}
{Hermann, Paris, 1959.}

\bib{shtsvl}{I.E. Shparlinski, M.A. Tsfasman, S.G. Vladut,}
{Curves with Many Points and Multiplication in Finite Fields.}
{Lectures Notes in Mathematics, Springer-Verlag, 1518 (1992), 145-169.}

\bib{stich}{H. Stichtenoth,}
{Algebraic Function Fields and Codes.}
{Springer Universitext, Berlin/Heidelberg/New York, 1993.}

\bib{weil}{A. Weil,}
{The Field of Definition of a Variety.}
{American Journal Of Mathematics, 78 (1956), 509-524.}


\end{thebibliography}
\end{document}